\documentclass[12pt]{article}

\usepackage{amscd,amsmath, amssymb, fancyhdr, epsfig,color}
\definecolor{dblue}{rgb}{0,0,.6}
\usepackage[colorlinks=true, linkcolor=dblue, citecolor=dblue, filecolor = dblue, menucolor = dblue, urlcolor = dblue]{hyperref}

\numberwithin{equation}{section}

\newcommand{\version}{version 1.1,\ \   Apr. 18, 2019}

\def\eqref#1{(\ref{#1})}
\newcommand{\goth}{\mathfrak}
\newcommand{\g}{{\mathfrak g}}
\newcommand{\arrow}{{\:\longrightarrow\:}}
\newcommand{\Z}{{\Bbb Z}}
\newcommand{\C}{{\Bbb C}}

\newcommand{\R}{{\Bbb R}}

\def\1{\sqrt{-1}\:}
\newcommand{\restrict}[1]{{\left|_{{\phantom{|}\!\!}_{#1}}\right.}}
\newcommand{\cntrct}                % contraction with a vector field
{\hspace{2pt}\raisebox{1pt}{\text{$\lrcorner$}}\hspace{2pt}}

\makeatletter
\def\x@arrow{\DOTSB\Relbar}
\def\xlongequalsignfill@{\arrowfill@\x@arrow\Relbar\x@arrow}
\newcommand{\xlongequal}[2][]{%
        \ext@arrow 0099\xlongequalsignfill@{#1}{#2}}
\def\xlongrightarrowfill@{\arrowfill@\relbar\relbar\longrightarrow}
\newcommand{\xlongrightarrow}[2][]{%
        \ext@arrow 0099\xlongrightarrowfill@{#1}{#2}}
\makeatother

\renewcommand{\bar}{\overline}
\renewcommand{\phi}{\varphi}
\renewcommand{\epsilon}{\varepsilon}

\renewcommand{\leq}{\leqslant}

\newcommand{\End}{\operatorname{End}}
\newcommand{\Mat}{\operatorname{Mat}}

\newcommand{\Id}{\operatorname{Id}}
\newcommand{\Cl}{\operatorname{{\mathcal C}\!\ell}}

\newcommand{\Tr}{\operatorname{Tr}}

\newcommand{\Spin}{\operatorname{Spin}}

\newcommand{\CC}{\mathcal{C}}
\newcommand{\HH}{\mathcal{H}}
\newcommand{\bbZ}{\mathbb{Z}}
\newcommand{\bbQ}{\mathbb{Q}}
\newcommand{\bbR}{\mathbb{R}}
\newcommand{\bbC}{\mathbb{C}}

\newcommand{\bbP}{\mathbb{P}}

\newcommand{\emrp}{\mathrm{End}}

\newcommand{\Pin}{\mathrm{Pin}}

\newcommand{\sdot}{{\raisebox{0.16ex}{$\scriptscriptstyle\bullet$}}}
\newcommand{\frg}{\mathfrak{g}}

\newcommand{\frsl}{\mathfrak{sl}}

\newcommand{\frso}{\mathfrak{so}}
\newcommand{\gtot}{\mathfrak{g}_{\mathrm{tot}}}
\newcommand{\hrarr}{\hookrightarrow}
\newcommand{\wdg}{\wedge}
\newcommand{\Lef}{\mathrm{L}}
\newcommand{\cnv}{\,\lrcorner\,\,}
\newcommand{\fldk}{\mathsf{k}}

\newcounter{Mycounter}[section]
\newcounter{lemma}[section]
\renewcommand{\thelemma}{{Lemma \thesection.\arabic{lemma}}}
\newcommand{\lemma}{%
    \setcounter{lemma}{\value{Mycounter}}
    \refstepcounter{lemma}
    \stepcounter{Mycounter}
    {\noindent \bf \thelemma:\ }}

\newcounter{claim}[section]
\newcounter{sublemma}[section]
\newcounter{corollary}[section]
\renewcommand{\thecorollary}{{Corollary \thesection.\arabic{corollary}}}
\newcommand{\corollary}{%
    \setcounter{corollary}{\value{Mycounter}}
    \refstepcounter{corollary}
    \stepcounter{Mycounter}
    {\noindent \bf \thecorollary:\ }}

\newcounter{theorem}[section]
\renewcommand{\thetheorem}{{Theorem \thesection.\arabic{theorem}}}
\newcommand{\theorem}{%
    \setcounter{theorem}{\value{Mycounter}}
    \refstepcounter{theorem}
    \stepcounter{Mycounter}
    {\noindent \bf \thetheorem:\ }}

\newcounter{conjecture}[section]
\newcounter{proposition}[section]
\renewcommand{\theproposition}
      {{Proposition \thesection.\arabic{proposition}}}
\newcommand{\proposition}{%
    \setcounter{proposition}{\value{Mycounter}}
    \refstepcounter{proposition}
    \stepcounter{Mycounter}
    {\noindent \bf \theproposition:\ }}

\newcounter{definition}[section]
\renewcommand{\thedefinition}
      {{Definition~\thesection.\arabic{definition}}}
\newcommand{\definition}{%
    \setcounter{definition}{\value{Mycounter}}
    \refstepcounter{definition}
    \stepcounter{Mycounter}
    {\noindent \bf \thedefinition:\ }}

\newcounter{example}[section]
\newcounter{remark}[section]
\renewcommand{\theremark}{{Remark \thesection.\arabic{remark}}}
\newcommand{\remark}{%
    \setcounter{remark}{\value{Mycounter}}
    \refstepcounter{remark}
    \stepcounter{Mycounter}
    {\noindent \bf \theremark:\ }}

\newcounter{problem}[section]
\newcounter{question}[section]
\newcommand{\proof}{{\bf Proof:\:}}

\makeatletter

\@addtoreset{equation}{section} \@addtoreset{footnote}{section} \makeatother

\def\blacksquare{\hbox{\vrule width 5pt height 5pt depth 0pt}}
\def\endproof{\blacksquare}

\begin{document}
%%%%%%%%%%%%%%%%%%%%%%%%%%%%%%%%%%%%%%%%%%%%%%%%%%%%%%%%%%%%
\begin{center}
{\LARGE\bf
Kuga-Satake construction and
cohomology of hyperk\"ahler manifolds}
\medskip
%%%%%%%%%%%%%%%%%%%%%%%%%%%%%%%%%%%%%%%%%%%%%%%%%%%%%%%%%%%%

Nikon Kurnosov\footnote{Partially supported by the Russian Academic Excellence Project '5-100' within AG Laboratory NRU-HSE and by the contest ``Young Russian Mathematics''.},
Andrey Soldatenkov\footnote{Partially supported by
the SFB/TR 45 ``Periods, Moduli Spaces and Arithmetic of Algebraic
Varieties'' of the DFG (German Research Foundation).},
Misha Verbitsky\footnote{Partially supported by the Russian Academic Excellence Project '5-100' within AG Laboratory NRU-HSE,
FAPERJ E-26/202.912/2018 and CNPq - Process 313608/2017-2}

\end{center}

%%%%%%%%%%%%%%%%%%%%%%%%%%%%%%%%%%%%%%%%%%%%%%%%
{\small \hspace{0.02\linewidth}
\begin{minipage}[t]{0.85\linewidth}
{\bf Abstract}
%Given a compact K\"ahler manifold $M$,
%its cohomology $H^\sdot(M,\bbC)$ is equipped with
%an action of the Lie algebra $\gtot(M)$ ge\-nerated by
%all the Lefschetz $\frsl(2)$-triples. When $M$ is
%a compact hyperk\"ahler manifold of maximal holonomy,
%$\gtot(M)$ induces the Hodge decomposition for any
%complex structure of hyperk\"ahler type on $M$.
%We generalize the Kuga-Satake construction
%and produce for some $l\ge 0$ an embedding
%$\Psi\colon H^\sdot(M,\bbC)\hrarr H^{\sdot + l}(T, \bbC)$
%into the cohomology of a compact complex torus $T$ in such a way that $\gtot(M)$ 
%is embedded into $\gtot(T)$, and the embedding $\Psi$ % of $H^\sdot(M,\bbC)$ into $H^\sdot(T,\bbC)$
%is a morphism of $\gtot(M)$-modules
%(and hence a morphism of Hodge structures).
Let $M$ be a simple hyperk\"ahler manifold.
Kuga-Satake construction gives an embedding of $H^2(M,\bbC)$
into the second cohomology of a torus, compatible with the Hodge structure.
We construct a torus $T$ and an embedding
of the graded cohomology space $H^\sdot(M,\bbC) \to H^{\sdot+l}(T,\bbC)$
for some $l$, which is compatible
with the Hodge structures and the Poincar\'e pairing.
Moreover, this embedding is compatible with
an action of the Lie algebra generated by all
Lefschetz $\frsl(2)$-triples on $M$.
\end{minipage}
}
%%%%%%%%%%%%%%%%%%%%%%%%%%%%%%%%%%%%%%%%%%%%%%%%

{\small 
\tableofcontents
}

%%%%%%%%%%%%%%%%%%%%%%%%%%%%%%%%%%%%%%%%%%%%%%%%

\section{Introduction}

%%%%%%%%%%%%%%%%%%%%%%%%%%%%%%%%%%%%%%%%%%%%%%%%

%%%%%%%%%%%%%%%%%%%%%%%%%%%%%%%%%%%%%%%%%%%%%%%%%%%%%%%%%%%%%%%%%%%%
\subsection{Kuga-Satake construction and $k$-symplectic geometry}
%%%%%%%%%%%%%%%%%%%%%%%%%%%%%%%%%%%%%%%%%%%%%%%%%%%%%%%%%%%%%%%%%%%%%%

The classical Kuga-Satake construction \cite{_Kuga_Satake_}, introduced by
Ichiro Satake and Michio Kuga in 1967, is used to construct
a holomorphic embedding of the space of Hodge structures of K3 type
(identified with the bounded Hermitian symmetric domain of
Cartan  type IV)
into the space of Hodge structures of abelian varieties,
identified with the bounded Hermitian symmetric domain of
Cartan  type III. The construction attaches to a Hodge
structure $H$ of K3 type a Hodge structure $V$ of weight one,
where $V$ is the even part of the Clifford algebra of $H$ with
the intersection form.
This way, the Kuga-Satake construction produces
an abelian variety from a given polarized
K3 surface. More generally, from a non-algebraic K3 surface
the construction produces a compact complex torus.
%The Kuga-Satake construction was used by Deligne \cite{_Deligne_Weil_K3_}
%to prove Weil conjectures for K3 surfaces.

We give an interpretation of the
Kuga-Satake construction based on the theory
of $k$-symplectic spaces. Our main result
is \ref{_embe_compa_with_g-action_Theorem_}.

The $k$-symplectic spaces were introduced in
\cite{_SV:k-symplectic_}; however the corresponding
geometric structure for $k=3$ was known since
1990-ies under the name ``hypersymplectic''.
A hypersymplectic space (\cite{_Dancer_Swann:hypersy_})
is a vector space $V$ over a field $\fldk=\R$ or $\C$ equipped with a triple
of symplectic forms $\omega_1, \omega_2, \omega_3$
in such a way that the operators
$\omega_i \circ \omega_j^{-1}:\; V \arrow V$
generate the matrix algebra $\Mat(2, \fldk)$.

Hypersymplectic manifolds are manifolds equipped with
a triple of symplectic forms giving it a hypersymplectic
structure at each tangent space. They were introduced
by N. Hitchin \cite{_Hitchin:hypersymple_}, because
the hyperk\"ahler reduction can be naturally extended
to the hypersymplectic case. In fact, the hypersymplectic
manifolds admit a unique torsion-free connection
preserving the triple of symplectic forms (see
\cite{_JV:Instantons_}). As shown in
\cite{_JV:Instantons_},
any hyperk\"ahler manifold $M$ admits a natural
complexification $M_\C$, which is realized as a component
in the moduli of rational curves in its twistor space,
and this space is a complex hypersymplectic manifold.

In the paper \cite{_Kamada:hypersymplectic_},  H. Kamada
classified hypersymplectic structures on compact surfaces,
and gave interesting examples of hypersymplectic
structures on a Kodaira surface.\footnote{Kodaira
surface is a non-K\"ahler compact complex surface obtained as 
an isotrivial holomorphic elliptic fibration over an 
elliptic curve.}

In \cite{_Jardim_Verbitsky:trisymple_}, Hitchin's
construction for hyperk\"ahler and hypersymplectic
reduction was applied to a complex hypersymplectic
manifold obtained as a component in the space
of rational curves in the twistor space of
a hyperk\"ahler manifold. This procedure was 
called {\bf trihyperk\"ahler reduction}.

In \cite{_Jardim_Verbitsky:trisymple_}, trihyperk\"ahler
reduction was applied to the space of geometric instantons
on $\C P^3$. Using this method it was shown that the
moduli of geometric instantons in $\C P^3$ is smooth;  this solves a
longstanding conjecture by W. Barth.

In hyperk\"ahler and real hypersymplectic case, the moment
map associated with the $G$-action takes 
values in $\g^* \otimes_\R \R^3$, where
$\g^*$ is the dual space to the Lie algebra of $G$.
For the trihyperk\"ahler reduction, the moment map
takes values in $\g^* \otimes_\R \R^7$, suggesting
that there are higher analogues of symplectic, trisymplectic
and trihyperk\"ahler structures for which the
geometric reduction decreases the dimension by
$2^k \dim \g$; this higher analogues should be associated
with stable bundles on $\C P^n$, $n>3$, in the same way as
the hypersymplectic structures are associated to 
the instantons on $\C P^3$. 

Hypersymplectic structures were generalized to 
$k$-symplectic structures in \cite{_SV:k-symplectic_}.
The starting point was again hyperk\"ahler geometry,
but an entirely different facet. Consider a compact complex torus
$T$ in a sufficiently general deformation of a 
complex manifold $M$ of hyperk\"ahler type. In \cite{_SV:k-symplectic_}
it was shown that the first homology $H_1(T,\R)$
is equipped with a $k$-symplectic structure, 
where $k=b_2(M)$. The intuition behind this construction
is very simple. Consider a general triple of
cohomology classes $w_1, w_2, w_3 \in H^2(M)$;
it is not hard to see that their restrictions
to $H^2(T)$, considered as 2-forms on $H_1(T)$,
form a hypersymplectic structure. 

This is the
intuition which underlies the notion of a $k$-symplectic
structure, and here is its formal definition.

\hfill

%%%%%%%%%%%%%%%%%%%%%%%%%%%%%%%%%%%%%%%%%%%%%%%%%%%%%%%%%
\definition\label{_k_symple_Definition_} 
Let $V$ be a vector space, and $\Omega\subset
\Lambda^2 V^*$ a subspace of dimension $k$.
We say that $\Omega$ is a {\bf $k$-symplectic structure}
if the following two properties are satisfied

\begin{description}
\item[(i)] For any non-zero $w\in \Omega$, the form $w$ has maximal
rank, or its rank is $\frac 1 2 \dim V$. This implies,
in particular, that the set of degenerate forms $\omega
\in \Omega$ is a quadric (see \cite{_SV:k-symplectic_};
this is an easy exercise in linear algebra).
\item[(ii)] The quadric $Q \subset {\Bbb P}(\Omega)$ 
consisting of all degenerate forms $\omega \in \Omega$
is non-degenerate.
\end{description}

It is still unknown what kind of geometric structure
is given by a $k$-symplectic analogue of a hypersymplectic
manifold, when we are given $k$ symplectic forms on a
manifold $M$ satisfying the $k$-symplectic properties
at each tangent space $T_x M$.

%%%%%%%%%%%%%%%%%%%%%%%%%%%%%%%%%%%%%%%%%%%%%%%%%%%

\subsection{Hyperk\"ahler manifolds}

%%%%%%%%%%%%%%%%%%%%%%%%%%%%%%%%%%%%%%%%%%%%%%%%%%%

For the convenience of the reader we will recall below the main
definitions related to hyperk\"ahler manifolds.
For more details see \cite{GHJ} and \cite{_Besse:Einst_Manifo_}.

\hfill

\definition (E. Calabi, \cite{_Calabi_}) Let $(M, g)$ be a Riemannian
manifold with three integrable complex structure operators
$I, J, K \in \emrp(TM)$, satisfying the quaternionic relations
$I^2=J^2=K^2=IJK=-\Id$.  Suppose that $g$ is K\"ahler with respect
to $I$, $J$ and $K$. Then $(M, I, J, K, g)$ is called hyperk\"ahler.

\hfill

\definition A holomorphically symplectic manifold 
is a complex manifold equipped with a non-degenerate, holomorphic
$(2,0)$-form.

\hfill

The metric $g$ on a hyperk\"ahler manifold $M$ is K\"ahler
with respect to the complex structures $I$, $J$ and $K$. Hence we have
the K\"ahler forms $\omega_I$, $\omega_J$, $\omega_K$, each
of them of type $(1,1)$ with respect to the corresponding complex
structure. One can check that the form 
$\Omega_I:= \omega_J+\1\omega_K$ 
is a holomorphic symplectic (2,0)-form on
$(M,I)$. This shows that every hyperk\"ahler manifold is holomorphically
symplectic. The partial converse to this statement follows
from Yau's solution of the Calabi conjecture (see \cite{_Yau:Calabi-Yau_}):

\hfill

\theorem Let $M$ be a compact holomorphically symplectic K\"ahler
manifold. Then $M$ admits a hyperk\"ahler metric which is
uniquely determined by the cohomology class of its 
K\"ahler form.

\hfill

In this paper we will consider only compact K\"ahler holomorphically symplectic
manifolds. There is the following structure theorem for such manifolds.

\hfill

\definition A compact hyperk\"ahler manifold $M$ is
{\bf of maximal holo\-nomy}, or {\bf simple},
or {\bf IHS}, if $\pi_1(M)=0$, $H^{2,0}(M)=\C$.

\hfill

\theorem (Bogomolov's decomposition: \cite{_Bogomolov:decompo_})
Any compact hyper\-k\"ahler manifold admits a finite \'etale covering
which is a product of a torus and se\-veral 
maximal holonomy hyperk\"ahler manifolds.
The maximal holonomy hyperk\"ahler components
of this decomposition are uniquely determined.

\hfill

Further on, all hyperk\"ahler manifolds
are tacitly assumed to be of maximal holonomy.
The following theorem is crucial for the study
of cohomology algebras of IHS manifolds.

\hfill

\theorem
(Fujiki, \cite{_Fujiki:HK_})
Let $M$ be an IHS manifold, $\dim_\bbC M=2n$.
Then there exist a non-degenerate primitive integer quadratic form $q$ on $H^2(M,\Z)$
and a rational constant $c$, such that for any $\eta\in H^2(M,\bbC)$
we have $\int_M \eta^{2n}=c q(\eta,\eta)^n$.

\hfill

%%%%%%%%%%%%%%%%%%%%%%%%%%%%%%%%%%%%%%%%%%%%%%%%%%%%%%%%%%%%%%%
\definition\label{_BBF_Definition_}
The quadratic form $q$ from the theorem above is called
{\bf the Bogomolov-Beauville-Fujiki (BBF) form}.

\hfill

The BBF form is determined up to a sign
by the relation from Fujiki's theorem. The sign may be determined
by the condition $q(\Omega,\bar \Omega) > 0$, where $\Omega$ is the holomorphic symplectic form.
It is known that the form $q$ has signature $(3,b_2(M)-3)$.
It is negative definite on $\omega$-primitive forms, and positive
definite on $\langle \Omega, \bar \Omega, \omega\rangle$,
where $\omega$ is a K\"ahler form.

%%%%%%%%%%%%%%%%%%%%%%%%%%%%%%%%%%%%%%%%%%%%%%%%%%%%%%%%%%%%%%%%%%%%

\section{Graded Frobenius algebras and $\goth{sl}(2)$-triples}\label{_Lef_Fro_Section_}

%%%%%%%%%%%%%%%%%%%%%%%%%%%%%%%%%%%%%%%%%%%%%%%%%%%%%%%%%%%%%%%%%%%%%%

To relate $k$-symplectic structures
and the Kuga-Satake construction, we use
the Lie algebra action generated by all 
Lefschetz triples on the cohomology of
a hyperk\"ahler manifold. This action was
discovered and studied in
\cite{_Looijenga+Lunts_,_Verbitsky:cohomo_,_Verbitsky:coho_announce_}.
Let us now recall the main elements of these works.

\subsection{Frobenius algebras}

%%%%%%%%%%%%%%%%%%%%%%%%%%%%%%%%%%%%%%%%%%%%%%%%%%%%%%%%%
\definition\label{_Frobenius_Definition_}
Let $A^\sdot=\bigoplus_{i=0}^r A^i$ be a graded-commutative algebra with $\dim A^r=1$.
Consider the $A^r$-valued form on $A^\sdot$ mapping $x, y\in A^\sdot$ to the
$A^r$-component of $xy$. The algebra $A^\sdot$ is called {\bf degree $r$ graded
Frobenius algebra} if this pairing is non-degenerate. 

\hfill

The basic example of a Frobenius algebra
is the cohomology algebra of a compact manifold. 

A {\bf Lefschetz triple}
in a Frobenius algebra $A =\bigoplus_{i=0}^{2n} A^i$ is a triple
of operators $\Lef_\eta, \Xi, \Lambda_\eta\in \emrp(A^\sdot)$
where $\eta \in A^2$ is a fixed element, 
$\Lef_\eta(x):= \eta x$, $\Xi\restrict{A^i} = i-n$ and
$\Lambda_\eta$ is an element such that
$\Lef_\eta, \Xi, \Lambda_\eta$ form an $\goth{sl}(2)$-triple.
It is easy to see that such $\Lambda_\eta$ is uniquely
determined by $\Xi$ and $\eta$ (this statement is sometimes
called ``Morozov's lemma'', and sometimes included
in the statement of Jacobson-Morozov theorem).
Existence of one Lefschetz triple is a non-trivial condition; however,
the space of $\eta\in A^2$ for which the Lefschetz triple exists is Zariski open.
%We call $\eta\in A^2$ for which the Lefschetz triple
%exists a {\bf Lefschetz-type element in $A^2$}.

\hfill

\definition
The {\bf Frobenius-Lefschetz algebra} is a Frobenius
algebra admitting a Lefschetz $\goth{sl}(2)$-triple.

\hfill

\remark
Let $A^\sdot$ be a Frobenius-Lefschetz algebra.
Consider the Lie algebra $\g(A^\sdot)\subset \emrp(A^\sdot)$
generated by all $\goth{sl}(2)$-triples. It was computed for
the cohomology of hyperk\"ahler manifolds
in \cite{_Verbitsky:cohomo_,_Looijenga+Lunts_} and for other 
geometric examples of Frobenius algebras (flag varieties, Hodge classes on an abelian
variety) in \cite{_Looijenga+Lunts_}.

\hfill

\theorem\label{_intro_g_constru_Theorem_}
Let $M$ be a hyperk\"ahler manifold of maximal holonomy, $A^\sdot$ its
cohomology algebra and $\gtot(M):=\g(A^\sdot)$ the Lie algebra 
generated by all Lefschetz $\goth{sl}(2)$-triples.
Then $\gtot(M)$ is isomorphic to $\goth{so}(4,b_2(M)-2)$.

\hfill

\proof See \cite{_Looijenga+Lunts_} or the proof of \ref{_g_constru_Theorem_}; an
explicit description of this algebra is also given there. \endproof

\hfill

An important property of the Lie algebra 
$\gtot(M)\cong \goth{so}(4,b_2(M) -2)$ acting on cohomology of
a hyperk\"ahler manifold is that it contains the
Weil operators $W_I$ for all complex structures
of hyperk\"ahler type on $M$ (see Subsection \ref{_supersy_hk_Subsection_}). 
Recall that for a compact K\"ahler manifold $(M,I)$
{\bf the Weil operator} $W_I\in \End(H^\sdot(M))$ acts on $H^{p,q}(M)$ as $\1(p-q)$.

This implies that the Hodge structure on $H^\sdot(M,\bbC)$ is induced
by the action of $\gtot(M)$. We will use this observation to
construct the Kuga-Satake-type embedding of the cohomology of
$M$ into the cohomology of a complex torus (see Section \ref{_Clifford_Section_}
and \ref{_embe_compa_with_g-action_Theorem_}).
The construction essentially consists of three steps:

\begin{enumerate}
\item Start from a hyperk\"ahler manifold of maximal holonomy,
and let $H = H^2(M, \C)$ be its second cohomology space equipped with the
BBF form. Construct a $k$-symplectic space $V$ such that $\phi\colon H \hrarr \Lambda^2 V^*$
is its $k$-symplectic structure. Both $V$ and $\phi$ will be naturally defined
over $\bbQ$. We will interpret $\Lambda^2 V^*$ as the second cohomology of
a complex torus $T$ well-defined up to isogeny.
\item Let $\g(H)\subset \gtot(T)$ be the subalgebra generated by the image of $H$
in $\Lambda^2V^*$. Then $\g(H)\simeq \gtot(M)$. Consider the action of $\g(H)$ on $H^\sdot(T,\bbC)$.
Let $\g_0(H)\subset \g(H)$
be the set of all elements which preserve the grading. By construction,
the embedding $\phi$ is $\g_0(H)$-invariant.
Since the Weil endomorphisms of $H^2(T,\bbC)$ and of $H$ 
belong to $\g(H)$, the embedding $\phi$ is compatible with the Hodge structure.
\item One can choose $V$ in such a way that the $\gtot(M)$-module $\Lambda^\sdot V^*$
contains $H^\sdot(M,\bbC)$ as a submodule. This gives the desired embedding of Hodge structures.
\end{enumerate}

\subsection{$F$-algebras}

%%%%%%%%%%%%%%%%%%%%%%%%%%%%%%%%%%%%%%%%%%%%%%%%%%%%%%%%%%%%%
Let $V$ be a vector space over an algebraically closed field of
characteristic zero. Let $q$ be a non-degenerate scalar product on $V$,
and $S^\sdot V$ the symmetric algebra. We will identify $V$ and $V^*$ via $q$.
Then $q$ can be considered an element of $S^2 V$.
Multiplication by $q^k$ gives a natural embedding
$S^l V \hrarr S^{l+2k} V$. Denote by 
$R_{n,k}(V)\subset S^{n+k} V$ the 
orthogonal complement to the image of $S^{n-k}V \stackrel{\cdot q^k}\hookrightarrow S^{n+k}V$
with respect to the non-degenerate symmetric pairing on $S^{n+k}V$
induced by $q$. Let $F_n^\sdot(V)$ be the quotient of $S^\sdot V$ 
by the ideal generated by $\bigcup_k R_{n,k}(V)$, with the grading
multiplied by two, so that $F_n^{2i}(V)$ is the quotient of $S^iV$.

\hfill

\definition\label{_F_algebra_Definition_}
The algebra $F_n^\sdot(V)$ is called {\bf the $n$-th $F$-algebra of $V$}.
This is an even-graded algebra, with $\dim F^{4n}(V)=1$.

\hfill

\remark
By definition, the natural map $S^iV \arrow F_n^{2i}(V)$ is
an isomorphism for $i\leq n$. For degree greater than $n$, 
one has $F_n^{2n+2i}(V)\cong S^{n-i}V$. It is clear from this
description that $F_n^{2n+2i}(V)$ is dual to $F_n^{2n-2i}(V)$,
that is, $F_n^\sdot(V)$ is a graded Frobenius algebra.

\hfill

\remark
It is easy to see (see e. g. \cite{Bog}) that the 
$F$-algebra has the following description:
$$ 
F_n^\sdot(V)\cong \frac{S^\sdot V}{\langle x^{n+1} \mid x\in V,\,
q(x,x) = 0\rangle}.
$$

\hfill

\theorem
Let $F_n^\sdot(V)$ be the $F$-algebra of a vector space $V$,
and $\g$ the Lie algebra generated by all $\goth{sl}(2)$-triples.
Then $\g$ is isomorphic to the Lie algebra
$\g=\goth{so}(\tilde V)$, where $\tilde V=V_0 \oplus V_2 \oplus V_4$
is the Mukai extension of $(V,q)$ defined as in the
Section \ref{_Mukai_extension_Subsection_}.
Moreover, $F_n^\sdot(V)$ is an irreducible $\g$-representation
generated by $1\in F_n^0(V)$.

\hfill

{\bf Proof:} This is just an explicit description of 
irreducible components of the symmetric power of the fundamental representation of 
$\g=\goth{so}(\tilde V)$; see \cite{_Verbitsky:coho_announce_}
for more details.
\endproof

%%%%%%%%%%%%%%%%%%%%%%%%%%%%%%%%%%%%%%%%%%%%%%%%

\section{Clifford modules and the Kuga-Satake
  construction}
\label{_Clifford_Section_}

%%%%%%%%%%%%%%%%%%%%%%%%%%%%%%%%%%%%%%%%%%%%%%%%

\subsection{Clifford algebras}

We start by fixing some notations involving Clifford
algebras. For more details see \cite{LM}. In this section
$\fldk$ will always denote a field of characteristic zero.

Let $H$ be a finite-dimensional vector space over $\fldk$
and $q\in S^2H^*$ a non-degenerate symmetric bilinear form.
Let $T^\sdot H$ denote the tensor algebra and $I\subset T^\sdot H$ be the two-sided ideal
generated by all elements of the form $v\otimes v-q(v,v)\cdot 1$ for $v\in H$.
The Clifford algebra is by definition $\Cl(H,q) = T^\sdot H/I$.
When $H$ and $q$ are clear from the context, we will denote the Clifford algebra by $\CC$.
Recall that $\CC$ is $\bbZ/2\bbZ$-graded $\CC=\CC^0\oplus\CC^1$, and denote by $\alpha\colon \CC\to \CC$
the parity involution. Denote by $\beta\colon \CC\to \CC^{\mathrm{op}}$ the
antiautomorphism induced by $v_1\otimes\ldots\otimes v_k\mapsto v_k\otimes\ldots\otimes v_1$.
We will use the notation $\bar{a}=\alpha\beta(a)$.

Recall that we have the canonical embedding $H\hrarr \CC$. The Clifford group
is by definition $G = \{a\in \CC^\times\,|\, \alpha(a)Ha^{-1} = H\}$.
The group $G$ comes with the natural action on $H$. Note that any non-isotropic $x\in H$ is contained in $G$.
The action of $x$ on $H$ is by reflection: $v\mapsto -xvx^{-1} = v - 2q(x,x)^{-1}q(x,v)x$.

For any $a\in G$ and $v\in H$
we have $q(\alpha(a)va^{-1},\alpha(a)va^{-1}) = (\alpha(a)va^{-1})^2=av^2a^{-1}=q(v,v)$, so
the action of $G$ on $H$ is orthogonal, and we get a morphism $\rho\colon G\to \mathrm{O}(H,q)$.
This morphism is surjective because $\mathrm{O}(H,q)$
is generated by reflections, and they are in the image of $\rho$.
One can check that the kernel of $\rho$ consists of invertible scalars $\fldk^\times$ (\cite{LM}, Proposition 2.4).

Let $G^i= G\cap \CC^i$.
The surjectivity of $\rho$ and the description of its kernel imply that any element $a\in G$
is a product $a = x_1\cdot\ldots\cdot x_n$ where $x_i\in H$ are non-isotropic vectors.
It follows that $G = G^0\coprod G^1$ and $G^0$ is a normal subgroup of index two.
Moreover, if we define $N(a) = a\bar{a}$ then
$N(a) = (-1)^n\prod_{i=1}^n q(x_i,x_i) \in \fldk^\times$ and $N(a)$ is called the norm of $a$.
The group $\Pin(H,q)$ is by definition the kernel of $N$.
We have $\Pin = \Pin^0\coprod \Pin^1$ and by definition $\Spin(H,q) = \Pin^0(H,q)$.

The Lie algebra $\frso(H,q)$ can be canonically identified with $\Lambda^2H$ and
embedded into $\CC$ via the map $\Lambda^2H\to\CC$, $x\wdg y\mapsto \frac{1}{2}(xy-yx)$.
This identifies $\frso(H,q)$ with the subspace of $\CC$ spanned by the commutators
of elements of $H$. The Lie bracket with any element of this subspace preserves $H$ and this
gives the description of the canonical action of $\frso(H,q)$ on $H$.

\subsection{Invariant forms on $\CC$-modules}

Consider a $\bbZ/2\bbZ$-graded $\CC$-module $V=V^0\oplus V^1$ with a bilinear form
$$b\colon V\otimes V\to k.$$

\definition
The bilinear form $b$ is called $\CC$-invariant if $b(au,v) = b(u,\bar{a}v)$ for all
$a\in \CC$ and $u,v\in V$.

\hfill

\remark
Our convention in this definition differs from what can be found in
the literature by the parity automorphism. It makes no difference when
we consider only the action of $\CC^0$ on $V$, but in what follows
it is important to use precisely the above definition. 
%In the case $\fldk=\bbC$ the existence and properties of $\CC$-invariant bilinear
%forms are discussed in \cite{Del}. Our convention differs from \cite{Del} by a parity
%automorphism. 

\hfill

If $b$ is $\CC$-invariant, then for any $a\in G$ we
have $b(au,av) = b(u,\bar{a}av) = N(\bar{a})b(u,v) = N(a)b(u,v)$. In particular for $a\in \Pin(H,q)$
we have $N(a)=1$, and we see that $b\in (V^*\otimes V^*)^{\Pin(H,q)}$.

We will now consider $\CC$ as a left module over itself. There exists a
natural $\CC$-invariant symmetric bilinear form on $\CC$. For $u,v\in \CC$ let $$\tau(u,v) = \Tr(u\bar{v}).$$
Here $\Tr(a)$ is the trace of $a$ considered as an endomorphism of $\CC$.
Note that the following relations hold: $\Tr(a) = \Tr(\alpha(a)) = \Tr(\beta(a)) = \Tr(\bar{a})$; $\Tr(uv) = \Tr(vu)$;
$\Tr(a)=0$ for $a\in\CC^1$.
This clearly implies that $\tau(u,v) = \tau(v,u)$. Moreover, for any $a,u,v\in \CC$ we have
$\tau(au,v) = \Tr(au\bar{v})=\Tr(u\bar{v}a) = \Tr(u\overline{\bar{a}v}) = \tau(u,\bar{a}v)$, so
the form $\tau$ is $\CC$-invariant. The trace of an odd degree element is zero, so the
form $\tau$ is even: $\tau\in S^2(\CC^0)^*\oplus S^2(\CC^1)^*$.

\hfill

\definition We will call the form $\tau$ described above the {\bf invariant trace form} on $\CC$.

\hfill

%%%%%%%%%%%%%%%%%%%%%%%%%%%%%%%%%%%%%%%%%%%%%%%%%%%%%%%%
\lemma\label{lem_trace_form}
Let $H$ be a quadratic vector space, $a\in H$ be a non-isotropic element,
$H'\subset H$ the orthogonal complement of $a$ and $\CC'=\Cl(H',q)$.
Let $\tau'$ be the corresponding trace form on $\CC'$.
Then the Clifford algebra decomposes into $\tau$-orthogonal direct sum
$\CC\simeq\CC'\oplus a\CC'$. We have $\tau|_{\CC'}=2\tau'$ and $\tau|_{a\CC'}=-2q(a,a)\tau'$.

\hfill

\proof
For $u,v\in \CC'$ we have $\tau(au,v) = \Tr(au\bar{v})$. The action of $au\bar{v}$
on $\CC'\oplus a\CC'$ exchanges the summands, so the trace is zero and the decomposition
$\CC\simeq \CC'\oplus a\CC'$ is $\tau$-orthogonal.

Assume that $u,v\in \CC'$ have the same parity. In this case $au\bar{v}=u\bar{v}a$ and the action of $u\bar{v}$ on $\CC'\oplus a\CC'$
is given by the matrix
$\left(\begin{array}{cc} u\bar{v} & 0 \\ 0 & u\bar{v} \end{array}\right)$.
We have $\tau(u,v) = \Tr(u\bar{v}) = 2\Tr'(u\bar{v})$ where $\Tr'$ is the trace in $\CC'$.
If $u$ and $v$ have different parity, the trace is zero.
This proves that $\tau|_{\CC'}=2\tau'$.

Note that $\tau(au,av) = \Tr(au\bar{v}\,\bar{a}) = -q(a,a)\Tr(u\bar{v})$.
This implies that $\tau|_{a\CC'}=-2q(a,a)\tau'$.
\endproof

\hfill

\corollary\label{lem_trform}
The invariant trace form $\tau$ is non-degenerate.

\hfill

\proof
This follows from the previous lemma and induction on $\dim H$. For $\dim H = 1$
we have $\CC\simeq \fldk[x]/(x^2-d)$ for some $d\in \fldk^\times$ and the trace form
is clearly non-degenerate.
\endproof

\subsection{The embedding}\label{subsec_embedding}
Let $V$ be a $\CC$-module of dimension $4n$ with a
$\CC$-invariant symmetric bilinear form $\tau\in S^2V^*$.
Assume that the form $\tau$ is non-degenerate.
For $x\in H$ and $u,v\in V$ let $\omega_x(u,v) = \tau(xu,v)$. Since
$\bar{x} = -x$, we have $\omega_x(u,v) = -\tau(u,xv) = -\omega_x(v,u)$ and so $\omega_x \in \Lambda^2V^*$. We get a map
\begin{equation}\label{eqn_embedding}
\varphi\colon H\to \Lambda^2V^*,\quad x\mapsto \omega_x.
\end{equation}

We will use the parity-twisted action of $\Pin(H,q)$ on $\Lambda^2V^*$:
for $\omega\in \Lambda^2V^*$, $a\in\Pin(H,q)$ and $u,v\in H$ let $(a\cdot\omega)(u,v) = \omega(a^{-1}u,\alpha(a)^{-1}v)$.

\hfill

\lemma\label{lem_emb} 
In the above setting $\varphi$ is a
morphism of $\Pin(H,q)$-modules. For any non-isotropic $x\in H$ the
form $\omega_x$ is non-degenerate.

\hfill

\proof
For any $a\in \Pin(H,q)$ we have $\bar{a} = a^{-1}$.
For any $x\in H$ and any $u,v\in V$ we obtain $\omega_{\alpha(a)xa^{-1}}(u,v)
=\tau(\alpha(a)xa^{-1}u,v) = \tau(xa^{-1}u,\alpha(\bar{a})v) = \tau(xa^{-1}u,\alpha(a)^{-1}v) = (a\cdot\omega_x)(u,v)$,
which proves the first claim.

The second claim follows from non-degeneracy of $\tau$: for any $u\in V$
one can find $v\in V$ with $\tau(u,v)\neq 0$. Then $\omega_x(u,xv)=
\tau(xu,xv) = -\tau(u,x^2v) = -q(x,x)\tau(u,v)$ which is non-zero
for non-isotropic $x$.
\endproof

\hfill

Denote by $W$ the image of $H$ in $\Lambda^2V^*$ under the map $\varphi$.
Consider the polynomial $p\in S^{2n}W^*$ given by $p(\omega) = \omega^{\wedge 2n}$
(here we identify $\Lambda^{4n}V^*$ with $\fldk$).
Analogously to \cite[Lemma 3.10]{_SV:k-symplectic_} we can prove that $W$ is a $k$-symplectic structure.

\hfill

\lemma\label{lem_emb2}
In the above setting $p$ is proportional to $q^{n}$. For any isotropic $x\in H$ the 
form $\omega_x$ has rank $2n$.

\hfill

\proof We may assume that the base field is algebraically closed.
The polynomial $p$ is the restriction to $W$ of $P\in S^{2n}(\Lambda^2V^*)$ defined in
the same way: $P(\omega) = \omega^{\wedge 2n}$.
But $P$ is $\Spin(H,q)$-invariant since the action of $\Spin(H,q)$ on $\Lambda^{4n}V^*$
is trivial ($\Lambda^{4n}V^*$ is one-dimensional). By \ref{lem_emb} $H$ is embedded
into $\Lambda^2V^*$ as a $\Spin(H,q)$-submodule,
so $p$ is invariant under the $\Spin(H,q)$-action on $H$. This action factors through $\mathrm{SO}(H,q)$,
so $p$ is invariant with respect to the special orthogonal group. Classical invariant
theory (see \cite{CP}, Theorem 5.6) implies that $p$ is proportional
to a power of the quadratic form $q$. We can choose the isomorphism $\Lambda^{4n}V^*\simeq \fldk$ so that $p = q^n$.

Let $x\in H$ be isotropic and $\dim(\ker\omega_x) = 2r$.
Pick a non-isotropic element $y\in H$ such that $q(x,y)\neq 0$.
Consider the polynomial
$\tilde{p}(t) = p(\omega_x+t\omega_y)$. We have $\tilde{p}(t) = q(x+ty)^n$,
so $\tilde{p}$ must have zero of order $n$ at $t=0$. But
\begin{eqnarray}
\tilde{p}(t)& =& (\omega_x+t\omega_y)^{\wedge2n}=\nonumber \\
&=&{\textstyle{2n\choose r}}t^r\omega_x^{\wedge(2n-r)}\wedge\omega_y^{\wedge r} + 
{\textstyle{2n\choose r+1}}t^{r+1}\omega_x^{\wedge(2n-r-1)}\wedge\omega_y^{\wedge(r+1)}+\ldots+t^{2n}\omega_y^{\wedge 2n},\nonumber
\end{eqnarray}
and $\omega_x^{\wedge(2n-r)}\wedge\omega_y^{\wedge r}\neq 0$. So the order of zero at $t=0$ is $r$, hence $r=n$.
\endproof

\hfill

Let $A$ be the subalgebra in $\Lambda^\sdot V^*$ generated by $W\subset \Lambda^2 V^*$. It is clear
that $A$ is a quotient of the symmetric algebra $S^\sdot H$. We would like
to describe the corresponding ideal in $S^\sdot H$.

\hfill

\proposition\label{prop_subring} Assume that $\fldk=\bar{\fldk}$.
Let $A^\sdot$ be the subalgebra in $\Lambda^\sdot V^*$ generated by $W$, where 
$V$ is a Clifford module over $\Cl(H,q)$. Then $A^\sdot$ is isomorphic to
the $F$-algebra $F_n^\sdot(H)$ (see \ref{_F_algebra_Definition_}),
where $n = \dim_\fldk V$.

\hfill

\proof
We will identify $H$ and $H^*$ via $q$. The algebra $A$ is a quotient of $S^\sdot H$, let $I$
be the corresponding ideal. Note that $I$ contains all
elements of the form $x^{n+1}$ for isotropic $x\in H$, as follows from \ref{lem_emb2}.

Recall the following facts from representation theory of the orthogonal group $\mathrm{O}(H,q)$
(see \cite{How}).
The form $q$ defines the trace maps $\mathrm{tr}\colon S^{k}H\to S^{k-2}H$
whose kernels are called spaces of harmonic polynomials, we will denote them by $\HH^k$.
By definition $\HH^1 = H$. It is known that $\HH^k$ are irreducible $\mathrm{O}(H,q)$-modules.
We also have the following decomposition of $S^{k}H$ into irreducible $\mathrm{O}(H,q)$-modules:
$$
S^kH= \bigoplus_{r\ge 0} q^r\HH^{k-2r},
$$
where by $q$ we mean the element of $S^2H$, corresponding to the form $q$
under the isomorphism $H\simeq H^*$.

Note that for any $x\in H$ with $q(x,x)=0$ we have $\mathrm{tr}(x^k)=0$, hence $x^k\in \HH^k$.
Since $\HH^k$ is irreducible, we see that it is spanned by the elements of the form $x^k$
with isotropic $x$. Also note that in the decomposition of the tensor product $\HH^k\otimes \HH^1$
into irreducible submodules we have two terms $\HH^{k+1}\oplus\HH^{k-1}$ plus some other
components, not isomorphic to $\HH^m$ (this can be checked using highest weight theory).
So the image of the multiplication map $\HH^k\otimes \HH^1\to S^{k+1}H$ is contained in
$\HH^{k+1}\oplus q\HH^{k-1}$. One can easily check that the image actually coincides with this
direct sum. More generally, for any $k\ge m$ the image of the multiplication map
$\HH^k\otimes\HH^m\to S^{k+m}H$ is $\HH^{k+m}\oplus q\HH^{k+m-2}\oplus\ldots\oplus q^m\HH^{k-m}$.

Going back to our algebra $A^\sdot$, note that the ideal $I$ contains $\HH^{n+1}$. Let us denote
by $J$ the ideal generated by $\HH^{n+1}$, then we have $J\subset I$ and $A$ is the quotient
of $B^\sdot=S^\sdot V/J$. For $k\le n$ we have $B^k=S^kH$. From
the description of the multiplicative structure on irreducible components $\HH^m$ it is easy
to see that $B^{n+i}\simeq q^i B^{n-i}$ for all $i\ge 0$. In particular $B^{2n}\simeq \fldk$.

We claim that multiplication maps $B^k\otimes B^{2n-k}\to B^{2n}\simeq \fldk$ give
non-degenerate pairing on $B^\sdot$. In order to see this, take an element $x\in q^r\HH^{k-2r}$ for
some $r$, so that $x=q^r\xi$ with $\xi\in \HH^{k-2r}$. Multiplication gives us
surjective map $\HH^{k-2r}\otimes\HH^{k-2r}\to \HH^{2k-4r}\oplus q\HH^{2k-4r-2}\oplus\ldots\oplus q^{k-2r} k$.
Hence there exists an element $\eta\in \HH^{k-2r}$, such that $\xi\eta=q^{k-2r}$. Then for
$y=q^{n-k+r}\eta\in q^{n-k+r}\HH^{k-2r}\subset B^{2n-k}$ we have $xy=q^n$. This shows that the pairing
is non-degenerate. In particular, all non-trivial ideals of the algebra $B^\sdot$ contain
the element $q^n$ which generates $B^{2n}$.

Note that for the algebra $A^\sdot$ the graded component $A^{2n}$ is non-zero, because $(\omega_x)^{\wedge 2n}\neq 0$
for a non-isotropic $x\in H$.
But if $A^\sdot$ were a quotient of $B^\sdot$ by some non-trivial ideal, this ideal would contain $B^{2n}\simeq \fldk$ and
$A^{2n}$ would be trivial. We conclude that $A$ is isomorphic to $B^\sdot$.
\endproof

%%%%%%%%%%%%%%%%%%%%%%%%%%%%%%%%%%%%%%%%%%%%%%%%%%%%%%%%%%%%%%
\subsection{Mukai extension} 
\label{_Mukai_extension_Subsection_}
%%%%%%%%%%%%%%%%%%%%%%%%%%%%%%%%%%%%%%%%%%%%%%%%%%%%%%%%%%%%%%

As above, let $(H,q)$ be a quadratic vector space. Let $\tilde{H} = \fldk \oplus H\oplus \fldk$
be the graded vector space with direct summands of degree $0$, $2$ and $4$. Define the quadratic form $\tilde{q}$ on $\tilde{H}$:
let $\tilde{q}((a,x,b), (a',x',b')) = q(x,x') - ab' - a'b$, so that degree $0$ and degree $4$ summands
make up a hyperbolic plane which is orthogonal to $H$, and the restriction of $\tilde{q}$ to $H$ is $q$.
We will refer to $(\tilde{H},\tilde{q})$ as Mukai extension of $(H,q)$. The Lie algebra $\frg = \frso(\tilde{H},\tilde{q})$
has the induced grading $\frg = \frg_{-2}\oplus \frg_0\oplus \frg_2$. Its components can be described as follows.

\hfill

\lemma\label{lem_mukai_ext}
We have $\frg_2\simeq \frg_{-2}\simeq H$ as vector spaces, $\frg_0\simeq \fldk \oplus \frso(H,q)$
as subalgebra (the first summand is the center of $\frg_0$). The action of $\frg_0$
on $\frg_{-2}$ and $\frg_2$ is via the standard representation of $\frso(H,q)$. The Lie bracket of two elements $x\in \frg_{-2}$
and $y\in \frg_2$  is given by $[x,y] = (q(x,y),x\wdg y)\in \frg_0$, where we use the
natural isomorphism $\frso(H,q)\simeq \Lambda^2H$.

\hfill

\proof Any element $x\in H$ defines two endomorphisms of $\tilde{H}$ of degree $2$ and $-2$:
one is given by $(a,z,b)\mapsto (0,ax,q(x,z))$, the other by $(a,z,b)\mapsto (q(x,z),bx,0)$. It is
straightforward to check that these endomorphisms are contained in $\frg$. It is also easy to
check the stated commutator relations. Hence we get two embeddings $H\hrarr \frg_{-2}$ and
$H\hrarr\frg_2$. The generator of the center of $\frg_0$ acts as $(a,z,b)\mapsto (-2a,0,2b)$.
The embedding $\frso(H,q)\hrarr \frg_0$ is the obvious one.
All these embeddings are isomorphisms for dimension reasons.
\endproof

\subsection{The Lie algebra action}
Recall that for a quadratic vector space $(H,q)$ we denote by $\CC$ the 
Clifford algebra $\Cl(H,q)$. For any $\CC$-module $V$ with
a non-degenerate $\CC$-invariant symmetric bilinear form $\tau\in S^2V^*$
we have constructed the embedding $\varphi\colon H\to\Lambda^2V^*$ (see \ref{eqn_embedding}).

We will show that the embedding $\varphi$
induces an action of the Mukai-extended Lie algebra $\frso(\tilde{H},\tilde{q})$ on $\Lambda^\sdot V^*$
compatible with the grading.

For $x\in H$ we will denote its image in $\Lambda^2V^*$ by $\omega_x$. Consider the
endomorphisms $\Lef_x, \Lambda_x\in \emrp(\Lambda^\sdot V^*)$ given by $\Lef_x\eta=\omega_x\wdg \eta$,
$\Lambda_x\eta = \omega_x\cnv\eta$ (to define the convolution with $\omega_x$ we identify $V$ and $V^*$
via $\tau$). The action of $a\in\CC$ on $V$ can be extended to a derivation of the algebra $\Lambda^\sdot V^*$
which we will denote by $\delta_a$.
Denote by $\Xi$ the endomorphism of $\Lambda^\sdot V^*$ that acts on the homogeneous
component of degree $k$ as multiplication by $k -
\frac{1}{2}\dim V$.

\hfill

\theorem\label{_g_for_k_symple_Corollary_}
The Lie subalgebra of $\emrp(\Lambda^\sdot V^*)$ generated by $\Lef_x$, $\Lambda_x$ for all $x\in H$
is isomorphic to the orthogonal Lie algebra 
$\frso(\tilde{H},\tilde{q})$ of the Mukai extension $\tilde H$.

\hfill

\proof Recall that $\frso(H,q)\simeq\Lambda^2H\hrarr\CC$, $x\wdg y\mapsto \frac{1}{2}(xy-yx)$ is an
isomorphism of $\frso(H,q)$ and a Lie subalgebra of $\CC$ spanned by commutators of elements
of $H$. Hence the operators $\delta_{\frac{1}{2}(xy-yx)}$ span a subalgebra isomorphic to $\frso(H,q)$.
The statement of the theorem now follows from \ref{lem_mukai_ext} and \ref{lem_commut}
by identifying $\frg_{2}$ with the subspace spanned by all $\Lef_x$ for $x\in H$, $\frg_{-2}$ with the subspace
spanned by all $\Lambda_x$ and the unit element of the center of $\frg_0$ with $\Xi$.
\endproof

\hfill

\remark In the case $\fldk = \bbR$ and $(H,q)$ the cohomology of a hyperk\"ahler
manifold with BBF form, one can prove the above theorem by checking
that the $\goth{so}(4, 1)$-relations hold for
any non-degenerate 3-dimensional subspace $H_3\subset H$.
%This subspace is associated with a hypersymplectic
%structure which is a complexification of a hyperk\"ahler structure,
%and satisfies the same relations.

\hfill

\lemma\label{lem_commut}
For any $x,y\in H$ we have
\begin{equation}\label{_commu_in_k-symple_Equation_}
[\Lambda_x,\Lef_y] = \delta_{\frac{1}{2}(xy-yx)} + q(x,y) \Xi.
\end{equation}

\proof Below we will give an explicit argument, based on computations.
However, these computations are possible to avoid as follows:
first, one checks \eqref{_commu_in_k-symple_Equation_} in
a 3-symplectic space (this is true and well known because
for 3-symplectic space the algebra generated by $L_{\omega_i}$
and $\Lambda_{\omega_j}$ is known and described explicitly in
Appendix (\ref{_so_4_1_Remark_}). Then, one can notice that 
it is possible to find a 3-symplectic structure on $V$ containing 
any generic $x, y \in H\subset \Lambda^2 V^*$, and hence the
relation \eqref{_commu_in_k-symple_Equation_} is true
on a Zariski dense set of $x, y$.

We can assume that $\fldk=\bar{\fldk}$. Since non-isotropic vectors form a 
Zariski-open subset of $H$, it is enough to prove the identity assuming that
$x$ is non-isotropic. So we fix a pair of elements $x,y\in H$ with $q(x,x)\neq 0$.
Note that for any $u\in V$ we have $\tau(xu,u)=0$, so we can find a basis
$V=\langle e_1,f_1,\ldots,e_{2n},f_{2n}\rangle$ where $f_i = xe_i$ and all elements
of the basis are pairwise orthogonal. We can also assume that $\tau(e_i,e_i) = 1$
and then $\tau(f_i,f_i) = -q(x,x)$. In the dual basis $\omega_x = -q(x,x) \sum e_i^*\wdg f_i^*$.
Using $\tau$ to identify $V$ and $V^*$ we can write $\omega_x = \sum e_i\wdg f_i$.

Let $\eta\in\Lambda^k V^*$ and $\xi_1,\ldots,\xi_k\in V$. We have 
\begin{eqnarray}
&&(\Lambda_x\Lef_y\eta)(\xi_1,\ldots,\xi_k) = \sum_{p=1}^{2n}(\omega_y\wdg\eta)(e_p,f_p,\xi_1,\ldots,\xi_k) = \nonumber\\
&&\sum_{p=1}^{2n}\tau(ye_p,f_p)\eta(\xi_1,\ldots,\xi_k) + \sum_{p=1}^{2n}\sum_{i=1}^k(-1)^i\tau(ye_p,\xi_i)\eta(f_p,\ldots,\check{\xi_i},\ldots) +\nonumber\\
&&\sum_{p=1}^{2n}\sum_{i=1}^k(-1)^{i-1}\tau(yf_p,\xi_i)\eta(e_p,\ldots,\check{\xi_i},\ldots) +\nonumber\\
&&\sum_{p=1}^{2n}\sum_{i<j}(-1)^{i+j-1}\tau(y\xi_i,\xi_j)\eta(e_p,f_p,\ldots,\check{\xi_i},\ldots,\check{\xi_j},\ldots).\nonumber
\end{eqnarray}
The last term in this formula is equal to $(\Lef_y\Lambda_x\eta)(\xi_1,\ldots,\xi_k)$.
The sum of the second and third terms is $(\delta_{xy}\eta)(\xi_1,\ldots,\xi_k)$, where we
use that
$$\xi_i = \sum_j(\tau(\xi_i,e_j)e_j - \frac{1}{q(x,x)}\tau(\xi_i,f_j)f_j).$$
The sum $\sum_p\tau(ye_p,f_p) = \sum_p\tau(ye_p,xe_p)$ can be rewritten in terms of the
trace of $xy$ acting on $V$: we have $\Tr(xy) = \sum_i(e_i^*(xye_i) + f_i^*(xyf_i)) =
\sum_i(\tau(e_i,xye_i) - q(x,x)^{-1}\tau(f_i,xyf_i)) = - 2\sum_i\tau(xe_i,ye_i)$.
We have obtained the following formula:
$$
[\Lambda_x,\Lef_y] = \delta_{xy} - \textstyle{\frac{1}{2}}\Tr(xy)\mathrm{Id}.
$$
Using the identity $xy = \frac{1}{2}(xy-yx) + q(x,y)$ we see that $\Tr(xy) = q(x,y) \dim V$,
and since $\delta_{q(x,y)}\eta = kq(x,y)\eta$ we have finished the proof.
\endproof

\hfill

\corollary
For any non-isotropic $x\in H$ the operators $\Lef_x$, $\Xi$ and $-q(x,x)^{-1}\Lambda_x$
form an $\frsl_2$-triple in $\emrp(\Lambda^\sdot V^*)$.

\hfill

\proof
The identity $[\Lef_x,-q(x,x)^{-1}\Lambda_x] = \Xi$ follows from the lemma.
The identities $[\Xi,\Lef_x] = 2\Lef_x$ and $[\Xi,-q(x,x)^{-1}\Lambda] = -2(-q(x,x)^{-1}\Lambda_x)$
follow from the fact that $\Lef_x$ is of degree $2$ and $\Lambda_x$ is of degree $-2$.
\endproof

\hfill

We obtain an action of the Lie algebra $\frso(\tilde{H},\tilde{q})$ on $\Lambda^\sdot V^*$
for any $\CC$-module $V$ with an invariant symmetric bilinear form $\tau$. Consider a
pair of such modules $V_1$, $V_2$ with the forms $\tau_1$, $\tau_2$ and let $\varphi_i\colon H\to \Lambda^2V_i^*$
be two embeddings constructed as above. On $V = V_1\oplus V_2$ we have the form $\tau_1\oplus \tau_2$
and it is clear by construction that $\varphi\colon H\to \Lambda^2V^*$ equals the composition
of $(\varphi_1, \varphi_2)\colon H\to \Lambda^2V_1^*\oplus \Lambda^2V_2^*$ and the natural
embedding $\Lambda^2V_1^*\oplus \Lambda^2V_2^*\hrarr \Lambda^2V^*$. This implies that the natural isomorphism
$\Lambda^\sdot (V_1\oplus V_2)^*\simeq \Lambda^\sdot V_1^*\otimes \Lambda^\sdot V_2^*$ is compatible
with the action of $\frso(\tilde{H},\tilde{q})$.

\hfill

\theorem\label{thm_all_rep}
For any quadratic vector space $(H,q)$ and any representation
$W$ of $\frg = \frso(\tilde{H},\tilde{q})$ there exists a $\Cl(H,q)$-module
$V$ with invariant symmetric bilinear form $\tau$, such that $\Lambda^\sdot V^*$
contains $W$ as a $\frg$-submodule.

\hfill

\proof
{\bf Step 1:} Assume additionally that the representation $W$ is irreducible.
We denote $\CC=\Cl(H,q)$. Let $V_1$ denote $\CC$ as
the left $\CC$-module with $\tau_1$ the invariant trace form. The action of $\frso(\tilde{H},\tilde{q})$
on $\Lambda^\sdot V_1^*$ induces an action of the group $\Spin(\tilde{H},\tilde{q})$.
Let us prove that this action is faithful. The construction behaves naturally under
the base field extensions, so we may assume that $\fldk=\bar{\fldk}$. The $\Spin$-group
is connected, and we have to check that its center acts non-trivially on $\Lambda^\sdot V_1^*$.

The structure of the center of $\Spin$ depends on $\dim(H)$ modulo $4$. If the
dimension is odd, then the only non-trivial element of the center is $-1\in \Spin(\tilde{H},\tilde{q})\subset \Cl(\tilde{H},\tilde{q})$.
This element is also contained in $\Spin(H,q)\subset \Spin(\tilde{H},\tilde{q})$ and
it clearly acts non-trivially on $V_1^*\subset \Lambda^\sdot V_1^*$.
In the case when the dimension of $H$ is even, we have to consider one more element
of the center. It can be described as follows. Let $\tilde{H}=\langle e_1,\ldots,e_{2d},f_1,f_2\rangle$
be the orthonormal basis, such that the elements $f_1$ and $f_2$ span the orthogonal
complement to $H\subset \tilde{H}$. Then the central element is
$g=e_1\ldots e_{2d} f_1 f_2$. Observe that $g$ is the product of
$g_1 = e_1\ldots e_{2d} \in \Spin(H,q)$ and $g_2=f_1 f_2\in \Spin(U)$ where
$U$ is the hyperbolic plane. Under the isomorphism $\Spin(U)\simeq k^\times$
the element $g_2$ corresponds to $\sqrt{-1}$. The action of $g_1$ on $\bbP(V_1^*)$
is non-trivial and the action of $g_2$ is trivial. Hence the action of $g$ on $V_1^*$
is non-trivial.

Since $\Lambda^\sdot V_1^*$ is a faithful $\Spin(\tilde{H},\tilde{q})$-module,
by \cite{DMOS}, Proposition 3.1 (p. 40), any irreducible representation of $\Spin(\tilde{H},\tilde{q})$ is contained
in $(\Lambda^\sdot V_1^*)^{\otimes N}$ for $N$ big enough. Then we can take $V=V_1^{\oplus N}$.

{\bf Step 2:} Consider the case when $W$ is not necessarily irreducible.
We will prove that $W$ can still be embedded into $(\Lambda^\sdot V_1^*)^{\otimes N}$ if $N$ is big enough.
We use induction on the number of irreducible components of $W$. Let $W\simeq W'\oplus W''$,
where $W'$ is irreducible.

Let us denote by $L$ the trivial one-dimensional representation of $\Spin(\tilde{H},\tilde{q})$.
There exists $k$, such that the representation $U = (\Lambda^\sdot V_1^*)^{\otimes k}$ contains $L$, so that
$U \simeq U'\oplus L$ for some representation $U'$. On the other hand, $U \simeq \Lambda^\sdot ((V_1^*)^{\oplus k})$
is still a faithful $\Spin(\tilde{H},\tilde{q})$-module by the argument from step 1.
This implies that $U'$ is faithful. It follows that for big $N'$ the representation $U^{\otimes N'}
\simeq (\Lambda^\sdot V_1^*)^{\otimes kN'}$ contains both the irreducible representation $W'$ and $L$.
By induction, we can assume that $U^{\otimes N''}$ contains both $W''$ and $L$, for some $N''$.
Then $U^{\otimes N'}\otimes U^{\otimes N''} \simeq (\Lambda^\sdot V_1^*)^{\otimes k(N'+N'')}$ contains
$W'\oplus W''\simeq W$ and $L$. This completes the induction step.
\endproof

%%%%%%%%%%%%%%%%%%%%%%%%%%%%%%%%%%%%%%%%%%%%%%%%%%

\section{Multidimensional Kuga-Satake construction}

%%%%%%%%%%%%%%%%%%%%%%%%%%%%%%%%%%%%%%%%%%%%%%%%%%%

\subsection{The embedding}

Let us apply the results of Section \ref{_Clifford_Section_}
to the cohomology of hyperk\"ahler manifolds. This will give
a multidimensional analogue of the Kuga-Satake construction.

\hfill

%%%%%%%%%%%%%%%%%%%%%%%%%%%%%%%%%%%%%%%%%%%%%%%%%%%%%%%%%%%%%%%%%%%%%%
\theorem\label{_embe_compa_with_g-action_Theorem_}
Let $M$ be a hyperk\"ahler manifold.
There exists an integer $l\ge 0$, a complex torus $T$,
an embedding $\gtot(M)\hrarr \gtot(T)$ of Lie algebras,
and an embedding 
$\Psi\colon H^\sdot(M,\bbC)\hrarr H^{\sdot + l}(T,\bbC)$
of $\gtot(M)$-modules.
For each complex structure $I$ of hyperk\"ahler type on $M$
there exists a complex structure on $T$ such that 
$\Psi$ is a morphism of Hodge structures. 

\hfill

{\bf Proof:} Let $H=H^2(M,\bbR)$ and $q$ be the BBF form (see \ref{_BBF_Definition_}).
We know that $\gtot(M) = \goth{so}(\tilde H,\tilde q)$ (see \ref{_g_constru_Theorem_}).
We apply the construction from Section \ref{subsec_embedding} and \ref{thm_all_rep}
to the pair $(H,q)$ and to $W = H^\sdot(M,\bbR)$. We obtain a $\Cl(H,q)$-module $V$,
an embedding $H\hookrightarrow \Lambda^2 V^*$ (a $k$-symplectic structure on $V$)
and an embedding $W\hrarr \Lambda^\sdot V^*$ of $\gtot(M)$-modules.
Note that the grading on $W$ and $\Lambda^\sdot V^*$ is induced by the action
of $\gtot(M)$, so we have a degree $l$ morphism $\Psi\colon H^\sdot(M,\bbR)\hrarr \Lambda^{\sdot+l}V^*$
of graded vector spaces, for some $l$ (the shift of degree is due to the difference of
dimensions of $M$ and $V$). We define $T$ to be the quotient of $V$ by a lattice.
Note that both $V$ and the morphism $\Psi$ are defined over $\bbQ$, since the BBF form
is defined over $\bbQ$, and the constructions from Section \ref{_Clifford_Section_}
behave naturally under base change.

The Hodge decomposition on $H^\sdot(M,\bbC)$ is given by the
$\goth u(1)$-action which belongs to the Lie algebra
$\goth{so}(H,q)\subset \goth{so}(\tilde H,\tilde{q})$ (see Section \ref{_supersy_hk_Subsection_}).
The corresponding  element of $\goth{so}(H,q)$ is a skew-symmetric
matrix $\mu$ of rank two.
Let $P=\langle e_1, e_2\rangle$, where $e_1 = \mathrm{Re}\,\Theta$, $e_2 = \mathrm{Im}\,\Theta$
and $\Theta$ is a generator of $H^{2,0}(M)$, normalized so that $q(e_i,e_i) = 1$.
Then $\mu = e_1e_2 \in \frso(H,q)\subset \Cl(H,q)$. It acts trivially on the orthogonal
complement to the 2-dimensional subspace $P$, and as $\begin{pmatrix} 0 & 2 \\ -2 & 0\end{pmatrix}$ on $P$.
However, on each Clifford module $\mu$ acts as a complex structure, because $\mu^2=-1$ in the
Clifford algebra. Therefore, $\mu$ acts as a complex structure on $H^1(T,\bbR)$.

The weight decomposition under the $e^{t\mu}$-action on 
$H^\sdot(M,\bbC)$ and $H^\sdot(T,\bbC)$ gives the Hodge decomposition on these cohomology spaces,
hence $\Psi$ is compatible with the Hodge structures.
\endproof

\subsection{Polarization}

Like in the classical Kuga-Satake construction, in the case when
the hyperk\"ahler manifold $M$ admits a polarization, the complex torus
$T$ from \ref{_embe_compa_with_g-action_Theorem_} admits a polarization too.
The proof of this is essentially the same as in the classical case.

Let $H_\bbQ = H^2(M,\bbQ)$, $H = H_\bbQ\otimes \bbR$, $h\in H_\bbQ$ an ample class
and $H'\subset H$ the orthogonal complement of $h$ with respect to the BBF form $q$.
Recall from the proof of \ref{_embe_compa_with_g-action_Theorem_}
that $T$ is a quotient of $V = V_1^{\oplus N}$ for sufficiently big $N$,
where $V_1\simeq \Cl(H,q)$. To prove that $T$ is polarized it
is enough to construct a rational Hermitian form $\sigma\in \Lambda^2 V_1^*$.
Let $\CC'=\Cl(H',q|_{H'})$. Note that $\Cl(H,q)\simeq \CC'\oplus h\CC'$
and this induces the decomposition $V_1\simeq \CC'^{\oplus 2}$ of rational Hodge structures,
because the Hodge structure on $V_1$ is defined by
a Weil operator which lies in $\CC'$. Hence it is enough to produce a polarization
on $\CC'$.

We have reduced the problem to the following: $(H_\bbQ,q)$ is a rational
quadratic vector space, $q$ has signature $(2,n)$ for some $n$ and we need to construct
a polarization on $V = \Cl(H,q)$ where $H=H_\bbQ\otimes \bbR$. Let $e_1,e_2\in H_\bbQ$
span a two-dimensional subspace in $H$ on which $q$ is positive.
Assume that $q(e_1,e_2) = 0$ and denote $a = e_1e_2\in \Cl(H_\bbQ,q)$.
For $x,y\in \Cl(H,q)$ let $$\sigma_a(x,y) = \Tr(xa\bar{y}).$$

\hfill

\proposition Either $\sigma_a$ or $-\sigma_a$ defines a polarization on $V$.

\hfill

{\bf Proof:} It is easy to check that $\sigma_a$ is skew-symmetric,
non-degenerate and $\Pin(H,q)$-invariant. The complex structure on
$V$ is given by an element $\mu=e_1'e_2'\in \frso(H,q)\subset \Cl(H,q)$, where
$e_i'\in H$ are such that $q(e_i',e_i') = 1$, $q(e_1',e_2')=0$.
We need to check that $\sigma_a(x,\mu x)\neq 0$ for $x\neq 0$, because this would imply that
$\sigma_a$ is sign-definite.

Note that we can find an element $g\in \Spin^+(H,q)$ such that
$\mu$ is proportional to $gag^{-1}$. Here $\Spin^+(H,q)$ is the connected
component of the identity in the $\Spin$-group. The orbit of $\sigma_a$ under
the action of $\Spin^+(H,q)$ is connected, and all forms in this orbit are
non-degenerate. Hence the form $\sigma_a$
is sign-definite if and only if $\sigma_\mu$ is. So it is enough to prove
that $\sigma_\mu(x,\mu x)\neq 0$ for $x\neq 0$. 

We prove this by induction on $\dim H$. Consider the orthogonal
decomposition $H = P\oplus H'$, where $P=\langle e_1',e_2'\rangle$.
Pick a non-zero element $z\in H'$ and let $H''$ be the orthogonal
complement to $z$. Let $\CC''=\Cl(H'',q|_{H''})$, then $V\simeq \CC''\oplus z\CC''$.
One can check like in the proof of \ref{lem_trace_form} that $\sigma_\mu|_{\CC''} = 2\sigma_\mu''$
and $\sigma_\mu|_{z\CC''} = -2q(z,z)\sigma_\mu''$, where $\sigma_\mu''$ is the corresponding
form on $\CC''$. Since $q(z,z)<0$, it is enough to prove that $\sigma_\mu''$ is positive-definite,
and so we can replace $H$ by $H''$. This completes the induction step. It remains
to consider the case when $H=P$. In this case $\Cl(H,q)$ is the algebra of $2\times 2$ real matrices.
One can check that $\sigma_\mu$ is positive definite on this algebra by a straightforward computation.
\endproof

\hfill

\remark The form $\sigma_a$ which we use to define the polarization is not
contained in the image of the embedding $\varphi$ from \ref{eqn_embedding}.

\hfill

\appendix

%%%%%%%%%%%%%%%%%%%%%%%%%%%%%%%%%%%%%%%%%%%%%%%%%%%%%%%%%%%%%

\section{Supersymmetry in hyperk\"ahler geometry}
\label{_Supersymm_Section_}

%%%%%%%%%%%%%%%%%%%%%%%%%%%%%%%%%%%%%%%%%%%%%%%%%%%%%%%%%%%%%

We give a summary of results which are
used to study the cohomology of a hyperk\"ahler manifold.
We follow \cite{_V:Mirror_}
and \cite{_Verbitsky:HKT_}; see also \cite{_FKS_}.

%%%%%%%%%%%%%%%%%%%%%%%%%%%%%%%%%%%%%%%%%%%%%%%%%%%%%%%%%%%%%
\subsection{Supersymmetry on K\"ahler manifolds}
%%%%%%%%%%%%%%%%%%%%%%%%%%%%%%%%%%%%%%%%%%%%%%%%%%%%%%%%%%%%%

Let $(M, I, g)$ be a compact K\"ahler manifold, $\omega$ its K\"ahler form,
and $\Lambda^\sdot(M)$ its de~Rham algebra.
On $\Lambda^\sdot(M)$, the following operators are defined:

\begin{enumerate}
\item de~Rham differential $d$, its adjoint $d^*$ and the Laplacian $\Delta$;

\item the Lefschetz operator $\Lef(\alpha)= \omega\wedge\alpha$ and  its adjoint $\Lambda(\alpha) := *\,\Lef * \alpha$;

\item the Weil operator $W\restrict{\Lambda^{p,q}(M)}=\1(p-q)$.
\end{enumerate}

\hfill

The following theorem is a convenient way to summarize
the K\"ahler relations and the Lefschetz $\goth{sl}(2)$-action.

\hfill

%%%%%%%%%%%%%%%%%%%%%%%%%%%%%%%%%%%%%%%%%%%%%%%%%%%%%%%%%
\theorem\label{_kah_susy_Theorem_}
These operators generate a Lie superalgebra
$\goth a$ of dimension $(5|4)$, 
acting on $\Lambda^\sdot(M)$. Moreover, the Laplacian $\Delta$ is
central in $\goth a$, hence $\goth a$ also acts on the
cohomology of $M$. 
\endproof

\hfill

\corollary
The ${\goth {sl}}(2)$-action of $\langle 
\Lef, \Lambda, \Xi\rangle$ and the action of the Weil operator $W$
commute with Laplacian, hence preserve the harmonic
forms on a K\"ahler manifold.

%%%%%%%%%%%%%%%%%%%%%%%%%%%%%%%%%%%%%%%%%%%%%%%%%%%%%%%%%%%%%%
\subsection{Supersymmetry on hyperk\"ahler manifolds}
\label{_supersy_hk_Subsection_}
%%%%%%%%%%%%%%%%%%%%%%%%%%%%%%%%%%%%%%%%%%%%%%%%%%%%%%%%%%%%%%

Let $(M, I, J,K, g)$ be a hyperk\"ahler manifold, $\omega_I$,
$\omega_J$, $\omega_K$ its K\"ahler forms.
On $\Lambda^\sdot(M)$, the following operators are defined:

\begin{enumerate}
\item de~Rham differential $d$, its adjoint $d^*$ and the Laplacian $\Delta$;

\item The Lefschetz operators
$$\Lef_I(\alpha)= \omega_I\wedge \alpha, \qquad \Lef_J(\alpha)= \omega_J\wedge \alpha, \qquad \Lef_K(\alpha)= \omega_K\wedge \alpha$$
and their adjoints
$$\Lambda_I(\alpha) = *\, \Lef_I * \alpha,\qquad \Lambda_J(\alpha)= *\, \Lef_J * \alpha, \qquad \Lambda_K(\alpha)= *\, \Lef_K * \alpha;$$

\item The Weil operators $W_{I}\restrict{\Lambda^{p,q}(M,I)}=\1(p-q)$, 
$W_{ J}\restrict{\Lambda^{p,q}(M,J)}=\1(p-q)$, $W_{K}\restrict{\Lambda^{p,q}(M,K)}=\1(p-q)$.
\end{enumerate}

The following result is an analogue of \ref{_kah_susy_Theorem_}:

\hfill

\theorem
These operators generate a Lie superalgebra
$\goth a$ of dimension $(11|8)$,
acting on $\Lambda^\sdot(M)$. Moreover, the Laplacian $\Delta$ is
central in $\goth a$, hence $\goth a$ also acts on the
cohomology of $M$. 
\endproof

\hfill

\remark The Weil operators span the Lie algebra
$\goth{su}(2)$ of unitary quaternions. 
In fact, one has $[\Lef_I, \Lambda_J]=W_K$,
 $[\Lef_J, \Lambda_K]=W_I$, $[\Lef_I, \Lambda_K]=-W_J$
(see \cite{_so(5)_}). This means that 
the quaternionic action belongs to $\goth a$.
The twisted de Rham differentials 
$d_I, d_J, d_K$, associated to $I,J,K$ also belong to
$\goth a$: $d_I= [ W_I, d]$, $d_J= [ W_J, d]$, 
$d_K= [ W_K, d]$.

\hfill

\remark
The action of $\goth a_{even}$ on $\goth a_{odd}$
is the fundamental representation of $\goth{sp}(1,1, {\Bbb H})$ in 
${\Bbb H^2}$, with the quaternionic Hermitian 
metric on $\goth a_{odd}$ provided
by the anticommutator. This implies that
the odd part \[ \langle d, d_I, d_J, d_K, d^*, d_I^*, d_J^*, d_K^*\rangle\]
generates the 9-dimensional odd Heisenberg algebra, with the
only non-trivial supercommutators being 
\[ \{d, d^*\}=\{d_I, d^*_I\}=\{d_J, d^*_J\}=\{d_K, d^*_K\}=\Delta.\]

\hfill

\remark The even part of $\goth a$
is isomorphic to $\goth{sp}(1,1, {\Bbb H})\oplus\R \cdot \Delta$.
The natural homomorphism to $\goth{sp}(1,1, {\Bbb H})\oplus\R \cdot \Delta$
is given by the previous remark. To see that it is an isomorphism 
we use the dimension count.

\hfill

%%%%%%%%%%%%%%%%%%%%%%%%%%%%%%%%%%%%%%%%%%%%%
\remark\label{_so_4_1_Remark_}
The algebra $\goth{sp}(1,1, {\Bbb H})$ is isomorphic
to $\goth{so}(4,1)$. To construct this isomorphism, we
consider the action of $\goth a$ on the 5-dimensional
graded vector space $H^\sdot=H^0\oplus H^2 \oplus H^4$, with
$H^0=\R$, $H^2=\langle \omega_I, \omega_J, \omega_K\rangle$,
$H^4=\R$, defined in the same way as the action of
$\goth a$ on the 5-dimensional subspace in the cohomology of a 
hyperk\"ahler manifold of real dimension 4.

\hfill

\remark
A Cartan subalgebra of \[ \goth{so}(4,1)\subset {\goth a}\]
can be given as $\langle \Xi, \1 W_I\rangle$.
The weight decomposition on $H^\sdot(M)$ 
associated with this Cartan algebra action 
coincides with the Hodge decomposition.

%%%%%%%%%%%%%%%%%%%%%%%%%%%%%%%%%%%%%%%%%%%%%%%%%%%%%%%%%%%%%%%%%
\subsection{Lie algebra generated by all $\goth{so}(4,1)$.}
%%%%%%%%%%%%%%%%%%%%%%%%%%%%%%%%%%%%%%%%%%%%%%%%%%%%%%%%%%%%%%%%%

The following was observed in \cite{_Verbitsky:cohomo_} (see also \cite{_Looijenga+Lunts_}).

\hfill

\proposition
Let $\Lef_\omega, \Xi, \Lambda_\omega$ and $\Lef_{\omega'},\Xi, \Lambda_{\omega'}$ 
be two $\goth{sl}(2)$-triples on a hyperk\"ahler manifolds.
Then $[\Lambda_{\omega'}, \Lambda_\omega]=0$.

\hfill

\proof From the Torelli theorem for hyperk\"ahler manifolds
(see e. g. \cite{_Americ-Verbitsky:Teich_s_}) it follows that the set of pairs
$\omega_I, \omega_J$ associated with the hyperk\"ahler
structures is Zariski dense in the space of all pairs
$\omega_1, \omega_2\in H^2(M,\R)$ such that
$q(\omega_1)= q(\omega_2)>0$ and $q(\omega_1, \omega_2)=0$.
Therefore it suffices to prove this relation for
$\omega, \omega'\in \langle \omega_I, \omega_J, \omega_K\rangle$.
In this case, the relation $[\Lambda_{\omega'}, \Lambda_\omega]=0$
follows from commutation relations in $\goth{so}(4,1)$.
\endproof

\hfill

From this result, the following structure 
theorem can be deduced (\cite{_Verbitsky:coho_announce_}, \cite{_Looijenga+Lunts_}).

\hfill

%%%%%%%%%%%%%%%%%%%%%%%%%%%%%%%%%%%%%%%%%%%%%%%%%%%%%%%%%%%%%%%%%
\theorem\label{_g_constru_Theorem_}
The algebra $\g$ generated by all $\goth{so}(4,1)$ 
for all hyperk\"ahler triples on a given
hyperk\"ahler manifold $M$ of maximal holonomy 
is isomorphic to $\goth{so}(4,b_2(M)-2)$.

\hfill

{\bf Proof. Step 1:}
Consider the action of $\g$ on the Mukai extension of $H^2(M,\bbR)$ (see Section \ref{_Mukai_extension_Subsection_}): 
\[ \tilde H:= \R\cdot x \oplus H^2(M,\bbR) \oplus \R\cdot y ,
\]
where $x$ has degree 0, $y$ has degree 4, $H^2(M,\bbR)$
is in degree $2$. On $\tilde H$ we have the extended quadratic form $\tilde{q}$.
The action of $\g$ is determined by the following properties:

\begin{enumerate}
\item It is compatible with the grading;

\item For all $\alpha, \beta\in H^2(M,\bbR)$, one has
$\Lef_\alpha x = \alpha$, $\Lef_\alpha \beta = q(\alpha, \beta)y$,
where $q$ is the BBF form.

\item $\Lambda_\alpha y = \alpha$, $\Lambda_\alpha \beta = q(\alpha, \beta)x$.
\end{enumerate}

To see that this action is well-defined, we need to check that
commutator relations hold. This follows from 
commutator relations in $\goth{so}(4,1)$ and Zariski density
of pairs $\alpha, \beta \in \langle \omega_I, \omega_J, \omega_K\rangle$
in the set of all pairs $\alpha, \beta\in H^2(M,\bbR)$. To obtain that 
the set of such pairs is Zariski dense, we use Torelli theorem
(see e. g. \cite{_Americ-Verbitsky:Teich_s_}). 

\hfill

{\bf Step 2:} We have constructed a homomorphism
$\g \arrow \End(\tilde{H})$. By construction, it preserves the
Mukai pairing $\tilde{q}\in S^2(\tilde{H}^*)$. This defines a homomorphism
$\Psi:\; \g \arrow \goth{so}(\tilde{H},\tilde{q})=\goth{so}(4,b_2(M) -2)$.

\hfill

{\bf Step 3:} $\Psi$ is surjective, because it is
surjective on generators, and injective, because the relations in
$\goth{so}(4,b_2 -2)$ can be obtained from relations in
$\goth{so}(4,1)$.
\endproof

\hfill

\remark 
It is easy to see that the center of the
corresponding Lie group $\Spin(4,b_2 -2)$ acts as $-1$  
on odd-dimensional cohomology and trivially on the 
even-dimensional cohomology (see \cite{_V:Mirror_}).

\hfill

{\bf Acknowledgement:} Part of this work was done by A.S. and M.V.
during the SCGP workshop on hyperk\"ahler manifolds in September 2016.
We are grateful to the organizers for the invitation and to the SCGP 
for perfect working conditions.

\hfill

{\small
\noindent
{\sc Nikon Kurnosov\\
{\sc Laboratory of Algebraic Geometry,\\
National Research University HSE,\\
Department of Mathematics, 
6 Usacheva Str. Moscow, Russia\\}
also: \\
{\sc
Department of mathematics, \\
University of Georgia,\\
1023 D. W. Brooks Drive,
Athens, GA 30602, USA\\}}

\noindent
{\sc Andrey Soldatenkov\\
{\sc Institut f\"ur Mathematik, Humboldt-Universit\"at zu Berlin,\\ Unter den Linden 6, 10099 Berlin\\
\verb|soldatea@hu-berlin.de|\\}

\noindent {\sc Misha Verbitsky\\
{\sc Instituto Nacional de Matem\'atica Pura e Aplicada\\
Estrada Dona Castorina, 110\\
Jardim Bot\^anico, CEP 22460-320\\
Rio de Janeiro, RJ - Brasil}\\ also: \\
{\sc
Laboratory of Algebraic Geometry, \\
National Research University HSE,\\
Department of Mathematics, 6 Usacheva Str. Moscow, Russia,}\\
\tt  verbit@impa.br},
}}


{\small 
\begin{thebibliography}{HKLR}


\bibitem[AV]{_Americ-Verbitsky:Teich_s_}
Amerik, E., Verbitsky, M., {\em Teichmuller space for hyperkahler
and symplectic structures}, preprint, arXiv:1503.01201, 2015.



%\bibitem[ABS]{ABS} Atiyah, M. F., Bott, R., Shapiro, A.,
%{\em Clifford modules},
%Topology, Vol. 3, 1964, suppl. 1, 3--38.

%\bibitem[Bea]{_Beauville_} 
% Beauville, A. {\em 
%Varietes K\"ahleriennes dont la premi\`ere classe de Chern est
%nulle.}  J. Diff. Geom. {\bf 18} (1983) 755 - 782.


\bibitem[Bes]{_Besse:Einst_Manifo_} 
Besse, 
A., {\em Einstein Manifolds}, Springer-Verlag, New York (1987)


\bibitem[Bo1]{_Bogomolov:decompo_}  
Bogomolov, F., {\em On the decomposition of 
K\"ahler manifolds with trivial canonical class}, Math. USSR-Sb.
{\bf 22} (1974) 580 - 583.


\bibitem[Bo2]{Bog}
Bogomolov, F.,
{\it On the cohomology ring of a simple hyper-K\"ahler manifold (on the results of {V}erbitsky)}
Geom. Funct. Anal., Vol. 6, No. 4, 1996, pp. 612--618.



\bibitem[C]{_Calabi_} 
Calabi,  E.,
{\em Metriques k\"ahleriennes et fibr\`es holomorphes}, 
Ann. Ecol. Norm. Sup. {\bf 12} (1979), 269-294.  

%\bibitem[De]{_Deligne_Weil_K3_}
%P. Deligne,
%{\em La conjecture de Weil pour les surfaces K3},
%Invent. Math., 15, 1972, pp. 206--226.

\bibitem[DMOS]{DMOS}
Deligne, P., Milne, J.S., Ogus, A., Shih, K.,
{\em Hodge cycles, motives, and Shimura varieties},
Lecture Notes in Mathematics, Vol. 900, Springer-Verlag, 1982.

\bibitem[CP]{CP}
de Concini, C., Procesi, C.,
{\em A Characteristic Free Approach to Invariant Theory}, Adv. in Math., 21, 1976, pp. 330--354.

\bibitem[DS]{_Dancer_Swann:hypersy_}
Dancer, A., Swann, A.,
Hypersymplectic manifolds.
In: {\em Recent developments in pseudo-Riemannian geometry}, 97--111,
ESI Lect. Math. Phys., Eur. Math. Soc., Z\"urich, 2008. 

\bibitem[FKS]{_FKS_}
Figueroa-O'Farrill, J. M.,
Koehl, C.. Spence, B.,
{\em Supersymmetry and the cohomology of (hyper)K\"ahler 
manifolds,} Nuclear Phys. B {\bf 503} (1997), no. 3, 614--626. 

\bibitem[F]{_Fujiki:HK_}  
Fujiki, A. {\em On the de Rham Cohomology Group of a Compact 
K\"ahler Symplectic Manifold}, Adv. Stud.
Pure Math. 10 (1987) 105-165.

\bibitem[GHJ]{GHJ} Gross, M., Huybrechts, D., Joyce, D.,
{\em Calabi-Yau manifolds and related geometries},
Universitext, Lectures from the Summer School held in Nordfjordeid, June 2001,
Springer-Verlag, Berlin, 2003, viii+239 pp.

\bibitem[Hi]{_Hitchin:hypersymple_}
Hitchin, N. J., {\em Hypersymplectic quotients}, 
Acta Acad. Sci. Tauriensis, supplemento al numero 124 (1990), 169--180.

\bibitem[How]{How}
Howe, R.,
{\it Remarks on classical invariant theory},
Trans. Amer. Math. Soc., Vol. 313, No. 2, 1989, pp. 539--570.

\bibitem[JV1]{_JV:Instantons_}
Jardim, M., Verbitsky, M.,
Moduli spaces of framed instanton bundles on $\C P^3$ and
twistor sections of moduli spaces of instantons on $\R^4$.
Adv. Math. {\bf 227} (2011),1526--1538.

\bibitem[JV2]{_Jardim_Verbitsky:trisymple_} 
Jardim, M., Verbitsky, M.,
{\em Trihyperk\"ahler reduction and instanton bundles on $\bbC P^3$},
arXiv:1103.4431 (to appear in Compositio Math.)

\bibitem[Ka]{_Kamada:hypersymplectic_}
Kamada, H., {\em Neutral hyper-K\"ahler structures on
  primary Kodaira surfaces}, 
Tsukuba J. Math.
      23 (2) (1999), 321-332. 100

%\bibitem[Ku]{_Kuga_}
% Kuga, Michio, {\em  Fiber varieties over a symmetric space whose
% fibers are abelian varieties} 1966 Algebraic Groups and
% Discontinuous Subgroups (Proc. Sympos. Pure Math.,
% Boulder, Colo., 1965) pp. 338-346 Amer. Math. Soc.,
% Providence, R.I. 

\bibitem[KS]{_Kuga_Satake_}
Kuga, M., Satake, I.,
{\em Abelian varieties attached to polarized K3-surfaces},
Math. Ann., 169 (1967), pp. 239--242.

\bibitem[LL]{_Looijenga+Lunts_}
Looijenga, E., Lunts, V., {\em A Lie algebra attached to a projective variety}, 
Invent. Math. 129 (1997), no. 2, 361-412.

\bibitem[LM]{LM}
Lawson, H.B., Michelsohn, M.-L.,
{\em Spin geometry}, Princeton University Press, 1989.

%\bibitem[MS]{_MS_} 
%A. Moroianu, U. Semmelmann,
%{\em Clifford structure on Riemannian manifolds},
%Adv. Math., 228, 2011, No. 2, 940--967.

%\bibitem[MP]{_Moroianu_Pilca_}
%Moroianu, Andrei; Pilca, Mihaela,
%{\em Higher rank homogeneous Clifford structures},
%J. Lond. Math. Soc. (2) 87 (2013), no. 2, 384-400. 

%\bibitem[Sa1]{_Satake_}
%Satake, I.
%{\em Holomorphic imbeddings of symmetric domains into a Siegel space}, 
%Amer. J. Math. 87 1965 425-461.

%\bibitem[Sa2]{_Satake_Clifford_}
%Satake, I.
%{\em Clifford algebras and families of abelian varieties,}
%Nagoya Math. J. 27 1966 435-446. 

\bibitem[SV2]{_SV:k-symplectic_}
Soldatenkov, A., Verbitsky, M., {\em $k$-symplectic
structures and absolutely trianalytic subvarieties in
hyperk\"ahler manifolds}, J. Geom. Phys. 92 (2015), 147-156.

\bibitem[V0]{_so(5)_} 
Verbitsky, M., 
{\em On the action of the Lie algebra 
$\frak{s}\frak{o}(5)$ on the cohomology
of a hyperk\"ahler manifold}, 
Func. Anal. and Appl. {\bf 24} (1990), 70-71.

\bibitem[V1]{_Verbitsky:cohomo_} 
Verbitsky, M., {\it Cohomology of 
compact hyperk\"ahler manifolds.}  alg-geom 
electronic preprint 9501001, 89 pages, LaTeX.

\bibitem[V2]{_Verbitsky:coho_announce_} 
Verbitsky, M.,
{\it Cohomology of compact hyperk\"ahler manifolds
and its applications,} GAFA vol. 6 (4) pp. 601-612 (1996).

\bibitem[V3]{_V:Mirror_}
Verbitsky, M.,
{\em Mirror Symmetry for hyperk\"ahler manifolds,}
 alg-geom/9512195,  Mirror symmetry, III (Montreal, PQ, 1995), 115--156,
   AMS/IP Stud. Adv. Math., 10, Amer. Math. Soc., Providence, RI, 1999.

\bibitem[V4]{_Verbitsky:HKT_}
Verbitsky, M., 
{\em 
Hyperk\"ahler manifolds with torsion, supersymmetry and Hodge theory},
Asian J. of Math., Vol. 6 (4), pp. 679-712, December 2002.

%\bibitem[V9]{_V:Torelli_}
%Verbitsky, M.,
%{\em A global Torelli theorem for hyperk\"ahler manifolds,}
% Duke Math. J. Volume 162, Number 15 (2013), 2929-2986.

\bibitem[Y]{_Yau:Calabi-Yau_} 
Yau, S. T., {\em On the Ricci curvature of a compact K\"ahler manifold 
and the complex Monge-Amp\`ere equation I.}  Comm. on Pure and Appl.
Math. 31, 339-411 (1978).


\end{thebibliography}
}
\end{document}